\documentclass[11pt]{amsart}

\usepackage[left=3cm, right=3cm, top=3cm, bottom=4cm]{geometry}
\usepackage{amsmath,amsfonts,amssymb,amsthm}
\usepackage{mathrsfs}
\usepackage{xcolor}
\usepackage{url}
\usepackage{graphicx}
\usepackage{subfigure}

\renewcommand{\phi}{\varphi}
\newcommand{\iy}{\infty}

\DeclareMathOperator{\trace}{Tr}

\DeclareMathOperator{\I}{I}

\newcommand{\Haar}{\textrm{Haar}}

\renewcommand{\leq}{\leqslant}
\renewcommand{\geq}{\geqslant}

\renewcommand{\d}[1]{\mathrm{d}#1}

\newcommand{\R}{\mathbb{R}}
\newcommand{\C}{\mathbb{C}}
\newcommand{\E}{\mathbb{E}}

\newcommand{\ol}[1]{\overline{#1}}

\newcommand{\scalar}[2]{\langle #1 , #2\rangle}
\newcommand{\crochet}[2]{\langle\!\langle \, #1 , #2 \, \rangle\!\rangle}

\newcommand{\braobsket}[3]{\langle #1 | #2 | #3 \rangle}

\renewcommand{\H}{\mathcal{H}}

\newcommand{\U}{\mathcal{U}}

\newcommand{\M}{\mathcal{M}}

\theoremstyle{plain}
\newtheorem{thm}{Theorem}[section]
\newtheorem{cor}[thm]{Corollary}
\newtheorem{prop}[thm]{Proposition}
\newtheorem{lem}[thm]{Lemma}

\theoremstyle{definition}
\newtheorem{defn}[thm]{Definition}

\theoremstyle{remark}
\newtheorem{rem}[thm]{Remark}

\title[Random pure states via UBM]{Random pure quantum states via unitary Brownian motion}
\author{Ion Nechita}
\address{CNRS, Laboratoire de Physique Th\'eorique, IRSAMC, Universit\'e de Toulouse, UPS, 31062 Toulouse, France.}
\email{nechita@irsamc.ups-tlse.fr}

\author{Cl\'ement Pellegrini}
\address{Institut de Math\'ematiques, IMT, Universit\'e de Toulouse (UMR 5219), 31062 Toulouse, Cedex 9, France.}
\email{clement.pellegrini@math.univ-toulouse.fr}

\keywords{}
\subjclass[2000]{}

\begin{document}

\begin{abstract}
We introduce a new family of probability distributions on the set of pure states of a finite dimensional quantum system. Without any a priori assumptions, the most natural measure on the set of pure state is the uniform (or Haar) measure. Our family of measures is indexed by a time parameter $t$ and interpolates between a deterministic measure ($t=0$) and the uniform measure ($t=\infty$). The measures are constructed using a Brownian motion on the unitary group $\mathcal U_N$. Remarkably, these measures have a $\mathcal U_{N-1}$ invariance, whereas the usual uniform measure has a $\mathcal U_N$ invariance. We compute several averages with respect to these measures using as a tool the Laplace transform of the coordinates.
\end{abstract}

\maketitle

\section{Introduction}

Defining models of randomness for  quantum objects has become a central problem in quantum information theory which has found many interesting and novel applications. Probability measures on the set of quantum states have been investigated thoroughly in recent years \cite{ZS01, Nec07, BZ06} both in the physical and the mathematical literature. Random quantum channels have also been a subject of interest \cite{HW08, CN10a, CN10b}. In particular, ensembles of quantum channels are the central idea behind the recent breakthroughs in the additivity conjecture \cite{Has09}. Random quantum states can arise in two different ways. First, they describe states of open systems which are subjected to random interaction with an (unknown) environment. This aspect has been the starting point of the so-called \emph{induced measures} which describe finite dimensional systems in interaction with a usually larger, but finite dimensional environment. Statistical ensembles of quantum states can also be used to study physical properties of generic states, such as entanglement, purity or other physically relevant quantities.

Defining a model of randomness for quantum states amounts to specifying a probability measure on the set of density matrices on the corresponding Hilbert space. When one considers only \emph{pure} states (rank one density matrices), it turns out that there exists a unique natural candidate for such a probability measure: the Lebesgue measure on the unit sphere of the underlying Hilbert space, $\pi_\infty$. This measure on the set unit vectors, called the \emph{Fubini-Study} or the \emph{uniform} distribution, is canonical in the sense that it is invariant under changes of bases: an element $ \psi$ of this ensemble has the same distribution as one of its rotations $U \psi$, for any given unitary matrix $U\in\mathcal U_N(\mathbb C)$. This invariance property, which characterizes the uniform distribution, justifies its use when no information about the internal structure of the system is known.

Very recently, new ensembles of pure states have been constructed, in order to take into account any available \emph{a priori} information on the system. In \cite{CNZ10}, the authors introduce an ensemble of states for quantum multipartite systems. Given a graph which encodes the entanglement between the different parties, a probability measure on the pure states of the total system is constructed. The ensemble of \emph{graph states} is different from the uniform ensemble because it contains the structure information about the initial entanglement present in the system.

In the present work, we are going to generalize the uniform ensemble in another direction, by removing the rotation invariance condition. Pure state ensembles on a single-partite system  with no \emph{full} invariance property exhibit \emph{preferred} states which have a larger probability than other states. The model introduced in the present work has a symmetry group of smaller dimension $\mathcal U_{N-1}(\mathbb C)$ and this feature makes is suitable for modeling systems on which some partial information is available. Given a fixed state $\psi$ of a quantum system, we shall introduce a one parameter family of probability measures  $\pi^\psi_t$ indexed by a real parameter $t \geq 0$. The parameter $t$ can be interpreted as a \emph{time parameter} is such a way that for $t=0$, the measure is deterministic, being supported  on the state $\psi$, and in the limit $t \to \iy$, the measure $\pi^\psi_t$ approaches the uniform measure $\pi_\infty$. For each value of the parameter $t$, the measure $\pi^\psi_t$ is invariant under the subgroup of rotations which leave invariant the vector $\psi$, making it the \emph{preferred} state of the measure. Our construction is based on the \emph{unitary Brownian motion}, a stochastic process valued in the set of unitary matrices. At fixed time $t$, this process itself is an interpolation between the identity matrix (at $t=0$) and the unique invariant Haar measure on the compact group of unitary matrices (when $t \to \iy$). The construction is motivated by the similar procedure that was used in the definition of the Fubini-Study measure.

The paper is organized as follows. In section \ref{sec:ubm} we review the definition and some basic properties of the unitary Brownian motion. Section \ref{sec:rps} contains the definition of the new family of ensembles of pure states, that are analyzed in section \ref{sec:PDE} using the Laplace transform. Finally, we compute in section \ref{sec:properties} averages of some quantities of interest in quantum information theory.
\bigskip

Let us now introduce some notation. In quantum information theory,  any norm one vector (or pure state) $x$ gives rise to a probability vector. More precisely, if $e_i$ is the canonical basis of $\H \simeq \C^N$, then we can decompose $\psi$ in the following form
$$\psi=\sum\psi_i e_i.$$
To such a state we naturally associate the probability vector
\begin{equation}\label{eq:prob-from-state}
	p(\psi) = (|\psi_1|^2, \ldots, |\psi_N|^2).
\end{equation}
Physically, if $e_i$ determines the level of energies of an atom and if $x$ represents the wave function describing this atom, the quantity $\vert \psi_i\vert^2$ represents the probability to be in the energy level $e_i$, that is
$$\mathbb{P}[\textrm{to\,\,be\,\,in\,\,the\,\,state\,\,}e_i]=\vert \psi_i\vert^2.$$

Another physical motivation related to probability theory concerns the measurement of observables. It is known in quantum mechanics, that a physical quantity of a quantum system $\H \simeq \C^N$ is represented by an observable, which is an auto-adjoint operator on $\H$. Let $A$ be an observable and $A=\sum_{i=1}^p\lambda_iP_i$ be its spectral decomposition. If $\psi$ is a reference vector state of $\H$, it follows from the axioms of quantum mechanics that a measurement of the observable $A$ gives a random result $\lambda_i$:
$$\mathbb{P}[\mathrm{to\,\, observe}\,\,\lambda_i]=\Vert P_i\psi\Vert^2.$$
In particular if the projectors $P_i$ are the one dimensional projectors on $\mathbb C e_i$ we recover the previous probability.

Let us now recall some elements of probability theory that we are going to use. If $X$ and $Y$ are two independent real Gaussian random variables of mean $0$ and variance $1/2$, then $Z = X + iY$ is said to have a complex Gaussian distribution of mean $0$ and variance $1$. We denote by $\mathcal N_\C(0,1)$ the law of $Z$. A complex vector $(Z_1, \ldots, Z_n)$ is said to have a multivariate complex Gaussian distribution $\mathcal N_\C^n(0,\I_n)$ if the random variables $Z_1, \ldots, Z_n$ are independent and have distribution $\mathcal N_\C(0,1)$.

We shall also extensively use the \emph{Haar} (or uniform) measure $\Haar_N$ on the unitary group $\U(\C)$; it is the unique probability measure which is invariant by left and right multiplication by unitary elements:
\begin{multline}\forall V, W \in \U_N(\C), \quad \forall f:\U_N(\C) \to \C \text{ Borel},\hfill\\ \quad \int_{\U_N(\C)}f(U) d\Haar_N(U) = \int_{\U_N(\C)}f(VUW) d\Haar_N(U).
\end{multline}

\section{Unitary Brownian Motion}\label{sec:ubm}

This section is devoted to the presentation of the unitary Brownian motion: definition, properties, stochastic calculus, invariant measure. In particular, we present all the ingredients that we are going to use for generating random quantum states.

The unitary Brownian motion refers to the natural definition of a Brownian motion on the unitary group of complex matrices. This is a special case of a Brownian motion on a differential manifold and more precisely on a compact Lie group (in differential geometry this is sometimes called the heat kernel measure).  

\subsection{Definition}

We now define the unitary Brownian motion (UBM). To start, let us introduce some notations. For $N \in \mathbb{N}$, we denote $\U_N(C)$ the unitary group on $ \M_N(\mathbb{C})$, that is, $\U_N(C)=\{U\in \textrm{GL}_n(\C)/UU^*=I\}$ and let denote by $\M_N^{sa}(\mathbb{C})$ the set of Hermitian matrices on $\M_N(\mathbb{C})$, that is, $\M_N^{sa}(\mathbb{C})=\{H\in \M_N(\mathbb{C})/H^*=H\}$. The set $\M_N^{sa}(\mathbb{C})$ is a real linear subspace of $\M_N(\mathbb{C})$, that we endow with the scalar product
$$\langle A,B\rangle=N\,\textrm{Tr}[A^*B]=N\,\textrm{Tr}[AB].$$
This way, we can consider the Brownian motion on $\M_N^{sa}(\mathbb{C})$, the unique Gaussian process $(H_t)$ which satisfies
\begin{equation}
\forall s,t,\,\forall A,B\in\M_N^{sa}(\mathbb{C}),\,\mathbb{E}[\langle A,H_s\rangle\langle B,H_t\rangle]=(s\wedge t)\langle A,B\rangle.
\end{equation}
In an equivalent way, the process $(H_t)$ has the same distribution of the random Hermitian matrix whose upper-diagonal coefficients are $\frac{1}{\sqrt{2N}}(B^{kl}_t+iC^{kl}_t)$ and whose diagonal coefficients are $\frac{1}{\sqrt{N}}D^k_t$, where $(B^{kl}_t,C^{kl}_t,D^{k}_t)$ are independent standard real Brownian motions.

We have now all the ingredients to define the unitary Brownian motion. This is the process $(U_t)_{t \geq 0}$, solution of the stochastic differential equation\footnote{This equation is written in Ito form.}
\begin{equation}
\left\{\begin{array}{ccc}dU_t&=&i(dH_t)\,U_t-\displaystyle{\frac{1}{2}}U_tdt\\
U_0&=&I.
\end{array}\right.
\end{equation}
In particular, we have
\begin{equation}
	d{U_t^*} = -i U_t^* d{H_t} - \frac{1}{2}U_t^* \d{t}.
\end{equation}
As a warm-up calculation, let us check that the process $(U_t)$ is unitary.  In order to compute $d(U_t^*U_t)$, we need to use Ito stochastic formulas for matrix valued stochastic processes. Such formulas have been derived in \cite[Section 2.1]{Ben11}: if $(X_t)$ and $(Y_t)$ are two matrix valued stochastic processes defined by
$$dX_t=A_tdH_tB_t+C_tdt,\quad dY_t=D_tdH_tE_t=F_tdt,$$
then we have the following Ito formula
\begin{equation}\label{eq:ito}
d(X_tY_t)=X_tdY_t+dX_t\,Y_t+\frac{1}{N}\textrm{Tr}[B_tD_t]A_tE_tdt.
\end{equation}
This way, we obtain
\begin{eqnarray} d\,U_tU_t^*&=&U_t^*d U_t+\frac{1}{n}\textrm{Tr}[U_tU_t^*]I\,d t\nonumber\\
&=& \left(\frac{1}{N}\textrm{Tr}[U_tU_t^*]I-U_tU_t^*\right)dt.
\end{eqnarray}
Since the unitarity condition is satisfied at $t=0$ ($U_0U_0^*=I$) and the identity is a solution of the above ordinary differential condition, by uniqueness of solutions we have $U_tU_t^*=I$ for all time $t$.

\subsection{Laplace Beltrami operator, Markov generator.}

For the sake of completeness, let us describe the Markov generator of this process. In particular, this allows us to motivate the definition of the Brownian motion from a geometric point of view. To this end, we denote by $\mathfrak{u}_N(\C)$ the Lie algebra of $\U_N(C)$  and by $(X_1,\ldots,X_{N^2})$ an orthonormal basis of ${\mathfrak{u}}(N)$. For all smooth functions $F:\U_N(C)\rightarrow \R$, we define
\begin{equation}
(\mathcal{L}_{X_i}F)(U)=\frac{d}{dt}_{\vert t=0}\left(F(Ue^{tX_i})\right),\quad\forall i.
\end{equation}

The operator $\frac{1}{2}\Delta$, where
\begin{equation}
\Delta=\sum_{i=1}^{n^2}\mathcal{L}^2_{X_i},
\end{equation}
is then the Markov generator of the unitary Brownian motion (this justifies the name heat kernel measure which is sometimes used to defining this process). The operator $\Delta$ is actually the Laplace-Beltrami operator on the Riemmanian manifold $\U_N(C)$ endowed with the Riemmanian metric induced by the scalar product on matrices $\langle A,B\rangle = N \trace[A^*B]$. Let us stress that this operator does not depend on any particular choice of an orthonormal basis.  The Markov generator character of $\Delta$ is expressed in the following proposition.
\begin{prop}  \cite[Proposition 2.1]{LM10} Let $F:\U_N(\C) \to \C$ be a function of class $C^2$. Then for all $t\geq0$, we have
\begin{eqnarray}
F(U_t)=F(I)+\sum_{i=1}^{N^2}\int_0^t(\mathcal{L}_{X_i}F)(U_s)d\langle X_k,iH_s\rangle+\int_0^t\frac{1}{2}\Delta F(U_s)ds
\end{eqnarray}
and the processes $(\langle X_k,iH_t\rangle),k=1,\ldots,N^2$ are independent standard real Brownian motions.
\end{prop}
The proof of this proposition relies on the classical Ito formula; this is a classical result of stochastic analysis on manifolds. A particular consequence of this proposition is that for all smooth function $F$, the process defined by
$$F(U_t)-F(U_0)-\int_0^t\frac{1}{2}\Delta F(U_s)ds,$$
for all $t$ is a martingale (with respect to the natural filtration associated with the Brownian motions $(B^{kl}_t,C^{kl}_t,D^{k}_t)$) and $U_t$ is the unique process "in distribution" satisfying such a property (this property characterizes the unitary Brownian motion and could be used as a starting definition).

\subsection{Invariant measure}

As announced in the introduction, the UBM will help us to define a new family of random states which interpolates between deterministic and uniformly distributed random states. This relies on the large time behavior and the invariant measure of the UBM. 

\begin{thm}\label{Haarinv}
For all initial unitary conditions $U_0$ the solution of the stochastic differential equation
\begin{equation}\label{eq:EDS-UBM}
dU_t=i(dH_t)\,U_t-\displaystyle{\frac{1}{2}}U_tdt,
\end{equation}
converges in distribution to the Haar measure $\Haar_N$ on $\U_N(\mathbb C)$. In other words, the Haar measure is the unique invariant measure of the Markov process, solution of \eqref{eq:EDS-UBM}.
\end{thm}

The theorem above justifies the property of interpolation between the identity operator ($U_0 =I$) and the Haar measure for large time ($t$ goes to infinity) for the unitary operator $U_t$. This property is essential for our definition of new ensembles of random pure states. This fact will be made precise in Section \ref{sec:rps}.

Theorem \ref{Haarinv} also shows that in general the distribution of the UBM $(U_t)$ is not invariant by unitary multiplication (except under the invariant measure) but the distribution of $(U_t)$ is nevertheless invariant by unitary conjugation and by inversion.

\begin{prop}\label{prop:UBM-invariance}
Let $(U_t)$ be the UBM defined by the SDE \eqref{eq:EDS-UBM} and let $V$ be any unitary matrix in $\U_N(\mathbb{C})$. The processes $(VU_tV^*)$ and $(U_t^{-1})$ have the same distribution than $(U_t)$.
\end{prop}

\begin{proof}
The property concerning the stochastic process $(VU_tV^*)$ follows from the fact that it satisfies the same stochastic differential equation \eqref{eq:EDS-UBM} with the same initial condition. Let $(W_t)$ defined by $W_t=VU_tV^*$ for all $t$, we have
\begin{eqnarray*}
dW_t&=&V(i\cdot dH_t\,U_t)V^*-\frac{1}{2}VU_tV^*dt\\
&=&i \cdot d(VH_tV^*)W_t-\frac{1}{2}W_tdt.
\end{eqnarray*}
Since $(VH_tV^*)$ is a Brownian motion on the Hermitian matrices, we see that $(W_t)$ and $(U_t)$ satisfies the same SDE. Hence, as they start with the same initial condition, the two processes must have the same distribution.

The statement for the inverse is an easy consequence of the Ito formula \eqref{eq:EDS-UBM} and the fact that the inversion corresponds to the complex adjoint and then is a linear mapping (not affected by derivation) \end{proof}

\subsection{Useful formulas}We continue by investigating some properties of the UBM which are going to be useful for studying the random pure states generated by the UBM. In particular, we will be interested in the properties of the coefficients of the matrix $(U_t)$ which we denote by $U_t^{jk}$, for $1 \leq j,k \leq n$. The stochastic differential equations satisfied by these elements are of the following form
\begin{equation}\label{eq:d-Ujk}
	d{U_t^{jk}} = i \sum_{s=1}^N (d{H_t^{js}})U_t^{sk} - \frac{1}{2}U_t^{jk} d{t}.
\end{equation}
and for their complex conjugates
\begin{equation}\label{eq:d-Ujk-bar}
	d{\ol{U_t^{jk}}} = -i \sum_{s=1}^N \ol{U_t^{sk}} d{H_t^{sj}} - \frac{1}{2} \ol{U_t^{jk}} d{t}.
\end{equation}
In the next section we will need the following expressions
\begin{align}\label{eq:d-modulus-j1}
	d\vert U_t^{j1}\vert^2&=d U_t^{j1}\,{\ol{U_t^{j1}}}+ U_t^{j1}d{\ol{U_t^{j1}}}+d U_t^{j1}d{\ol{U_t^{j1}}}\nonumber \\
	&=	i \sum_{s=1}^N (d{H_t^{js}})U_t^{s1}{\ol{U_t^{j1}}}-i \sum_{s=1}^N U_t^{j1}\ol{U_t^{s1}} d{H_t^{sj}}+\left(-\vert U_t^{j1}\vert^2+\frac{1}{N}\right)dt.		
	\end{align}
The previous formula relies on the stochastic bracket for the elements of $(H_t)$, that is,
\begin{equation}
	d{\crochet{H_t^{ij}}{H_t^{kl}}} = \frac{d t}{N}\delta_{il}\delta{jk},
\end{equation}
where $\delta$ is the Kronecker delta symbol.
From this, we obtain the brackets\footnote{Here we have adopted the notations $\crochet{}{}$ for the stochastic bracket not to be confused with the scalar product $\langle \, , \, \rangle$.} of the matrix coordinates which is given by
\begin{equation}\label{eq:d-bracket-modulus-j1-k1}
	d{\crochet{|U_t^{j1}|^2}{|U_t^{k1}|^2}} = \frac{2}{N} \left(|U_t^{j1}|^2 \delta_{jk} - |U_t^{j1}|^2|U_t^{k1}|^2\right)d{t}.
\end{equation}

\section{Random Pure States Generated by Unitary Brownian Motion}\label{sec:rps}

Before developing our theory for random pure states generated by UBM, we review the definition and the basic properties of the \emph{uniform} (or Fubini-Study) probability measure on the set of pure states (or unit vectors). The lack of any \emph{a priori} information on the state $\psi$ of a quantum system described by a Hilbert space $\H \simeq \C^N$ imposes the choice of a measure which should be  invariant by changes of bases. In our setting of finite dimensional complex Hilbert spaces, changes of bases are implemented by unitary operators $U \in \U_N(\mathbb C)$. As a consequence, we ask that the uniform probability measure should be unitarily invariant. A probability measure $\pi$ on the unit ball of $\H$ is called \emph{unitarily invariant} if for all Borel subsets $A$ and for all unitary operators $U \in \U_N(\mathbb C)$, 
\begin{equation}
\pi(UA) = \pi(A).
\end{equation}

The above condition determines uniquely the measure $\pi$: it is the normalized surface area of the unit ball of $\C^N$, which we shall denote by $\pi_\iy$, for consistency reasons which shall be clear later. Moreover, we introduce the image measure $\sigma_\iy = \mathrm{sq}_\# \pi_\iy$, where 
$\mathrm{sq} : \C^N \to \Delta_N$, $\mathrm{sq}[(\psi_i)] = (|\psi_i|^2)$. In other words, if the unit vector $\psi$ has distribution $\pi_\iy$, then the probability vector $(|\psi_i|^2)_{i=1}^N$ has distribution $\sigma_\iy$.

Other that the abstract definition of the invariant measure $\pi_\iy$, there are two more characterization of this probability that are important in what follows:
\begin{enumerate}
	\item Let $X \in \C^N$ be a standard complex Gaussian vector. Then $X / \|X\|$ has distribution $\pi_\iy$.
	\item Let $U \in \U_N(\mathbb C)$ be a Haar-distributed random unitary matrix. The first column (or any column, or any line) of $U$ has distribution $\pi_\iy$.
\end{enumerate}

The first statement above is useful when one needs to sample from $\pi_\iy$. The second statement above will be the starting point for the definition of new probability measures on the unit sphere of $\C^N$. 

Start with a fixed vector $\psi \in \C^N$ of norm one, and define the stochastic process $(\psi_t)$, where 
\begin{equation}
	\psi_t=U_t\psi,\quad \forall t \geq 0,
\end{equation}
and $(U_t)$ is a UBM starting at $U_0 = I_N$.
This gives rise to a stochastic process valued in the unit sphere, with $\psi_0 = \psi$. In the sequel, we study the properties of this process, whose distribution at time $t$ we denote by $\pi_t^\psi$. As before, the distribution of the probability vector $|\psi_t^j|^2$ is denoted by $\sigma_t^\psi$, making explicit the dependence in the initial condition $\psi_0 = \psi$. 

\begin{defn}\label{def:rps-ubm}
The ensemble of pure states (unit vectors) of $\C^N$ having distribution $\pi_t^\psi$ is called the \emph{unitary Brownian motion induced ensemble} at time $t$. The distribution of the square moduli of the coordinates of a random vector $\psi$ in this ensemble will be denoted by $\sigma_t^\psi$.
\end{defn}

Let us first discuss the connection between the distribution $\pi_t^\psi$ and the Haar measure $\pi_\iy$. First, note that the distribution of $\psi_t$ depends on the initial value $\psi$, in contrast with Haar distributed random states, whose distribution is invariant. In particular, we do not have for $\pi_t^\psi$ invariance by all unitary transformations, that is, if $V$ denotes a unitary operator, in general the processes $(V\psi_t)$ and $(\psi_t)$ have different distributions. More precisely we have the following result.

\begin{prop}\label{prop:invariance}
Consider a unitary operator $V$ and a (possibly random) unit vector $\psi$. Let $(\psi_t)$ be the process generated by a unitary Brownian motion independent of $\psi$, with initial condition $\psi_0=\psi$. Then, the processes $(V\psi_t)$ and $(U_t V\psi)$ have the same distribution. In particular, the processes $(V\psi_t)$ and $(\psi_t)$ have the same distribution if and only if $V\psi \sim \psi$. 
\end{prop}

\begin{proof}
This follows from the fact that
\begin{equation}\label{inv}
V\psi_t=U_t\psi = (VU_tV^*)V\psi = \tilde U_t V\psi .
\end{equation}
Since $(U_t)$ and $(\tilde U_t) = (VU_tV^*)$ have the same distribution and are independent of $\psi$ and $V\psi$, the first part is then straightforward. The second part is a trivial consequence of first part.
\end{proof}

\begin{cor}
Let $\psi$ be a random uniform unit vector, i.e. $\psi \sim \pi_\iy$. Consider an independent UBM $(U_t)$ which induces a process $(\psi_t)$ with initial condition $\psi$. Then, for all $t$, the vector $\psi_t$ has distribution $\pi_\iy$, i.e. $\pi_t^\psi=\pi_\iy$.
\end{cor}

Proposition \ref{prop:invariance} is more instructive when we look at deterministic initial conditions. In particular, the processes $( V\psi_t)$ and $(\psi_t)$ have the same distributions if and only if $ V\psi=\psi$, restricting the class of unitary transformations which leave invariant the process. In other words, the distribution of $( \psi_t)$ is invariant under all unitary transformations which fix the initial condition. In conclusion, the process $( \psi_t)$ has a $\U_{N-1}(\C)$ invariance group, whereas a uniform unit vector $\psi \sim \pi_\iy$ has a $\U_N(\C)$ invariance group.

\section{PDE for the Laplace transform of $\sigma_t^\psi$}\label{sec:PDE}

We are now in the position to develop our model in more detail, using as a main tool a partial differential equation satisfied by the Laplace transform of the amplitude vector. Let $\psi \in \C^N$ a fixed unit vector and consider the process generated by a UBM $(\tilde U_t)$ with $\tilde U_0 = I$:
\begin{equation}
	\psi_t = \tilde U_t \psi.	
\end{equation}
Consider also a fixed unitary matrix $V$ such that $V e_1 = \psi$, where $e_1 = (1, 0, \ldots, 0)$ is the first element of the canonical basis of $\C^N$. The one can write
\begin{equation}
	\psi_t = \tilde U_t V e_1 = U_t e_1,	
\end{equation}
where $(U_t)$ is another UBM starting at $U_0 = V$. Note that the choice of the matrix $V$ satisfying $V e_1 = \psi$ is not important, because of Proposition \ref{prop:invariance}. In this way, the initial condition of the problem has been transfered into the UBM $(U_t)$ and we have
$$\psi_t = U_t e_1=(U_t^{11},U_t^{21},\ldots,U_t^{N1}),$$
which corresponds to the first column of the unitary Brownian motion $U_t$. 
As is was discussed in the Introduction, such a unit norm vector gives rise to a probability vector $(|\psi_t^j|^2)_j = (|U^{j1}_t|^2)_j$ having distribution $\sigma_t^\psi$. This random variable is compactly supported, hence its Laplace transform determines its distribution. Let us define the Laplace transform by (notice the positive sign in the exponential)
\begin{equation}
	\phi(\lambda_1, \ldots, \lambda_N; t) = \E e^{\scalar{\lambda}{|U^{\cdot 1}_t|^2}} = \E \left[ \exp \sum_{j=1}^N \lambda_j |U_t^{j1}|^2 \right],
\end{equation}
for all $\lambda = (\lambda_1, \ldots, \lambda_N) \in\mathbb{C}^N$ and all $t\geq0$. Since the random variable is bounded, the function $\lambda \mapsto \phi(\lambda; t)$ is complex analytic for each $t$. 

The partial derivatives of the function $\phi$ read:
\begin{align}
\label{eq:d-j-phi}\partial_j \phi(\lambda; t) = \partial_{\lambda_j}\phi(\lambda; t) &= \E \left[ |U^{j1}_t|^2 \, e^{\scalar{\lambda}{|U^{\cdot 1}_t|^2}} \right];\\
\label{eq:d-jk-phi}\partial_{jk} \phi(\lambda; t) = \partial_{\lambda_j\lambda_k} \phi(\lambda; t) &= \E \left[ |U^{j1}_t|^2 |U^{k1}_t|^2 \, e^{\scalar{\lambda}{|U^{\cdot 1}_t|^2}}\right].
\end{align}

The following theorem is the main result of this paper, establishing a partial differential equation for the Laplace transform $\phi$. In principle, it allows to recover $\phi$ and then, by Laplace inversion, the probability vector $(|U_t^{j1}|^2)_{j=1}^N$.

\begin{thm}\label{thm:main}
The Laplace transform $\phi$ of the random vector $(|U^{j1}_t|^2),j=1,\ldots,N$ satisfies the following partial differential equation
\begin{equation}\label{eq:EDP-main}	\partial_t \phi = \frac{\sum_{j=1}^N \lambda_j}{N} \phi + \left<\frac{\lambda^2}{N} - \lambda,\nabla_\lambda \phi \right> - \frac{1}{N} \left<\lambda, H(\phi) \lambda\right>,
\end{equation}
where $\lambda^2$ represents the vector $(\lambda^2)_j =\lambda_j^2$ and $\nabla_\lambda$ and $H$ represent the gradient and the Hessian operators, i.e.
\begin{align*}
(\nabla_\lambda \phi)_j &= \partial_j \phi;\\
[H(\phi)]_{jk} &= \partial_j \partial_k \phi.
\end{align*}
\end{thm}
\begin{proof}
Using the multivariate Ito formula for the function $(x_1, \ldots, x_N) \mapsto e^{\scalar{\lambda}{x}}$ applied to the multidimensional process $(|U^{\cdot 1}_t|^2)$, we obtain:
\begin{equation}
	\d{\exp{\scalar{\lambda}{|U^{\cdot 1}_t|^2}}} = \sum_{j=1}^N \lambda_j  \exp{\scalar{\lambda}{|U^{\cdot 1}_t|^2}}\d{|U^{j1}_t|^2}   + \frac{1}{2}\sum_{j,k=1}^N \lambda_j\lambda_k \exp{\scalar{\lambda}{|U^{\cdot 1}_t|^2}} \d{\crochet{{|U^{j1}_t|^2}}{|U^{k1}_t|^2}}.
\end{equation}

Taking the expectation and using \eqref{eq:d-modulus-j1}, \eqref{eq:d-bracket-modulus-j1-k1}, we obtain
\begin{align}
	\d{\E \exp{\scalar{\lambda}{|U^{\cdot 1}_t|^2}}} = &\sum_{j=1}^N \lambda_j \E \left[ \exp{\scalar{\lambda}{|U^{\cdot 1}_t|^2}} \left(-\vert U_t^{j1}\vert^2+\frac{1}{N}\right)\right] dt \nonumber \\
	&+ \frac{1}{2}\sum_{j,k=1}^N \lambda_j\lambda_k \E \left[ \exp{\scalar{\lambda}{|U^{\cdot 1}_t|^2}} \frac{2}{N}\left(|U_t^{j1}|^2 \delta_{jk} - |U_t^{j1}|^2|U_t^{k1}|^2\right)\right] dt.
\end{align}

Using formulas \eqref{eq:d-j-phi}, \eqref{eq:d-jk-phi}, we obtain immediately the announced partial differential equation \eqref{eq:EDP-main} satisfied by $\phi$.
\end{proof}

From the above PDE one can obtain, in principle, all the information about the distribution $\sigma_t^\psi$ of the random vector $(|U^{j1}_t|^2)$. In the remainder of this section, we shall focus on the marginals $|U^{j1}_t|^2$ ($j$ fixed); covariances $|U^{j1}_t|^2|U^{k1}_t|^2$ and other statistical quantities will be investigated in the next section.

In the case of marginals, we compute in the following proposition the Laplace transform of the square modulus of one coordinate, in terms of the Kummer confluent hypergeometric function ${_1F_1}$, whose definition we recall :
\begin{equation}
{_1F_1}(a;b;z) = \sum_{k=0}^\infty \frac{(a)_k z^k}{(b)_k k!},
\end{equation}
where $(x)_n$ is the Pochhammer symbol, $(x)_n = x(x+1) \cdots (x+n-1)$.
\begin{prop}\label{prop-Lapl-transfo-j} Let $\phi_j(\lambda; t) = \phi(0, \ldots, 0, \lambda, 0, \ldots, 0; t) = \E\exp(\lambda|U^{j1}_t|^2)$ be the Laplace transform of the $j$-th coordinate of the first column of $U_t$. Then
\begin{equation}\label{eq:EDP-1-marginal}
	\partial_t \phi_j = \frac{\lambda}{N} \phi_j + \left( \frac{\lambda^2}{N} - \lambda \right) \partial_\lambda \phi_j - \frac{\lambda^2}{N} \partial_{\lambda\lambda} \phi_j,
\end{equation}
with the notation
\begin{equation}
	\partial_t \phi_j = \frac{\partial \phi_j}{\partial t}, \quad \partial_\lambda \phi_j = \frac{\partial \phi_j}{\partial \lambda},  \quad\partial_{\lambda\lambda} \phi_j = \frac{\partial^2 \phi_j}{\partial \lambda^2}.
\end{equation}
Given an initial condition $|U^{j1}_0|^2 = c \in [0, 1]$, there exists a sequence $(a_n)_{n \geq 0}$ of real numbers (depending on $c$) such that
\begin{equation}
	\phi_j(\lambda; t) = \sum_{n=0}^\iy a_n e^{-\Lambda_n t}  \lambda^n {_1F_1}(n+1; N+2n; \lambda),
\end{equation}
where
\begin{equation}
	\Lambda_n = n + \frac{n(n-1)}{N}.
\end{equation} 
\end{prop}
\begin{proof}
Equation \eqref{eq:EDP-1-marginal} follows from equation \eqref{eq:EDP-main} of Theorem \ref{thm:main} by letting $\lambda_k=0$ for all $k \neq j$. In the rest of the proof, we shall drop the index $j$, since the initial condition will be encoded into $\phi(\lambda;0) = \exp(\lambda|U^{j1}_t|^2)$.

Using separation of variables, we look for solutions of the form $\phi(\lambda; t) = f(\lambda)g(t)$. Neither of $f$ or $g$ can be zero, hence we obtain
\begin{equation}
	\frac{g'}{g} = -\frac{\lambda^2}{N}\frac{f''}{f} + \left(\frac{\lambda^2}{N} - \lambda \right) \frac{f'}{f} + \frac{\lambda}{N}.
\end{equation}
Note that the left hand side of the above equation depends only on $t$ and the right-hand since depends only on $\lambda$. This is impossible unless both are equal to a constant $C$, in which case $g(t) = e^{Ct}$ (we can move the constant factor to $f$) and $f$ satisfies the ordinary differential equation
\begin{equation}
	-\frac{\lambda^2}{N}f'' + \left(\frac{\lambda^2}{N} - \lambda \right) f' + \left(\frac{\lambda}{N}-C \right)f = 0.
\end{equation}
Writing $f$ as a power series $f(\lambda) = \sum_{n \geq 0} a_n \lambda^n$, we get:
\begin{align}
	C a_0 &= 0 \\
	a_1(1+C) &= \frac{1}{N} a_0 \\
	a_k\left( \frac{k(k-1)}{N} + k + C \right) &= \frac{k}{N} a_{k-1}  \qquad \forall k\geq 2.
\end{align}
These equations can be summarized as (we put $a_{-1} = 0$)
\begin{equation}
	a_k(C+\Lambda_k) = \frac{k}{N}a_{k-1} \qquad \forall k \geq 0,
\end{equation}

If $C \notin \{-\Lambda_k\}_{k \geq 0}$, it follows that $f=0$, which is impossible. Hence, $C = -\Lambda_n$ for some $n \geq 0$. We can compute all the coefficients of the series expansion of $f$ from the recurrence relations above:
\begin{align}
	a_m &= 0 \qquad 0 \leq m < n \\
	a_n &\text{ is free } \\
	a_{n+m} &=a_n \frac{(n+1)_m}{m!(N+2n)_m} = a_n \binom{n+m}{m}\frac{1}{(N+2n)_m} \qquad \forall m \geq 1.
\end{align}
We obtain the final expression for the Laplace transform:
\begin{equation}\label{eq:Laplace-phi}
	\phi(\lambda; t) = \E \exp(\lambda |U^{j1}_t|^2) = \sum_{n=0}^\iy a_n e^{-\Lambda_n t}  \sum_{m=0}^\iy \binom{n+m}{m}\frac{\lambda^{n+m}}{(N+2n)_m}.
\end{equation}
The second sum in the above formula admits a more compact expression using the Kummer confluent hypergeometric function $_1F_1$:
\begin{equation}
	\sum_{m=0}^\iy \binom{n+m}{m}\frac{\lambda^{n+m}}{(N+2n)_m} = \lambda^n {_1F_1}(n+1; N+2n; \lambda),
\end{equation}
and thus
\begin{equation}
	\phi(\lambda; t) = \sum_{n=0}^\iy a_n e^{-\Lambda_n t}  \lambda^n {_1F_1}(n+1; N+2n; \lambda).
\end{equation}
\end{proof}

\begin{rem}
Note that in the limit $t \to \iy$, only the $n=0$ term survives and we obtain
\begin{equation}
	\lim_{t \to \iy} \phi(\lambda; t) = \sum_{m=0}^\iy \frac{\lambda^m}{(N)_m}
\end{equation}
which is the result for the Haar measure. This is consistent with \cite{hpe}, Lemma 4.2.4, where the $m$-th moment in the Haar case was shown to be $\binom{N+m-1}{N-1}$.
\end{rem}

It is interesting to note that equation \eqref{eq:EDP-1-marginal} does not depend on the actual value of $j$. However, it does depend on the initial condition $c$. Next, we compute explicitly $\phi$ for particular values of the initial condition $c=|U_0^{j1}|^2\in[0,1]$. The values of the coefficients $a_n$ appearing in the proposition can be computed in principle from the following initial conditions:
\begin{align*}
\phi(\lambda; 0) &= e^{\lambda c}\\
\phi(0; t) &= 1.
\end{align*}
As a first observation, note that he latter relation fixes the value of the constant coefficient in the series, $a_0=1$. The first condition translates to ($p=n+m$):
\begin{equation}
	\sum_{n=0}^p a_n \binom{p}{n}\frac{1}{(N+2n)_{p-n}} = \frac{c^p}{p!}, \quad \forall p \geq 0.
\end{equation}
The above infinite triangular system of linear equations can be solved in principle and explicit formulas for the coefficients $a_n$ can be found. 

Analytical formulas can be obtained (and easily proved by induction) in two particular cases. For $c=0$, one can show that the unique solution to the equations above are given by
\begin{equation}
	a_n = \frac{(-1)^n}{(N+n-1)_n}, \quad \forall n \geq 0.
\end{equation}

Similarly, for $c=1$, one has
\begin{equation}
	a_n = \frac{(N-1)_n}{n! (N+n-1)_n}, \quad \forall n \geq 0.
\end{equation}

We were not able to obtain analytical expressions for all the coefficients in other particular cases. We gather next the first six coefficients in the important case $c=1/N$:
\begin{equation}
	a_0 = 1, \quad a_1 = 0, \quad a_2 = -\frac{-1+N}{2 N^2 (1+N)}, \quad a_3 = \frac{2 (-2+N) (-1+N)}{3 N^3 (2+N) (4+N)},
\end{equation}
\begin{equation}
	 a_4 = -\frac{(-1+N) \left(30-29 N+5 N^2\right)}{8 N^4 (3+N) (5+N) (6+N)}, \quad a_5=\frac{(-2+N) (-1+N) \left(84-79 N+7 N^2\right)}{15 N^5 (4+N) (6+N) (7+N) (8+N)}.
\end{equation}

\section{Properties of the measure $\sigma_t^\psi$}\label{sec:properties}

This section contains a list of results which address important statistical properties of probability vectors distributed along the measure $\sigma_t^\psi$. Moments, as well as covariances, are shown to satisfy ordinary differential equations that are solvable, see Propositions \ref{prop:moments} and \ref{prop:covariances}. We also compute quantities relevant to quantum information theory, such as average values of observables in Lemma \ref{lem:average-observable} and bounds for average R\'enyi entropies.

From Proposition \ref{prop-Lapl-transfo-j}, it is easy to obtain the expression of the moments of the $j$-th coordinate $|\psi^j_t|^2$. To this end, we define a family of maps $y_p:[0, \iy) \to [0, 1]$, $y_p(t) = \E |U_t^{j1}|^{2p}$. The maps, indexed by positive integers $p \geq 1$ depend implicitly on the size parameter $N$. The dependence on the index $j$ is encoded in the initial condition $y_p(0) = |\psi^j|^{2p} = |U_0^{j1}|^{2p}$. Interchanging the derivation operator $\partial^p$ and the expectation $\E$, we obtain
\begin{equation}
	y_p(t) = \frac{\partial^p \phi}{\partial\lambda^p}(0;t).
\end{equation}
The applications $y_p$ satisfies a particular system of ordinary differential equations.

\begin{prop}\label{prop:moments}
With the convention that $y_0 \equiv 1$, the applications $y_p$ satisfy the following ordinary differential equations on $\R_+$:
\begin{equation}\label{eq:EDP-moments}
	y_p' = -\Lambda_p y_p + \frac{p^2}{N}y_{p-1}
\end{equation}
The solution of the system \eqref{eq:EDP-moments} can be expressed in terms of $a_n$ and $\Lambda_n$ in the following way
\begin{equation}
	y_p(t) = \sum_{n=0}^p   \binom{p}{p-n}\frac{a_n }{(N+2n)_{p-n}}e^{-\Lambda_n t}.
\end{equation}
\end{prop}

\begin{proof}
Using Proposition \ref{prop-Lapl-transfo-j}, we can rewrite the formula \eqref{eq:Laplace-phi} for $\phi(\lambda; t) $ in the form
\begin{equation}
	\phi(\lambda; t) = \E e^{\lambda |U^{j1}_t|^2} = \sum_{p=0}^\iy \sum_{n=0}^p   \binom{p}{p-n}\frac{a_n e^{-\Lambda_n t}}{(N+2n)_{p-n}}\,\lambda^{p}
\end{equation}
and then, taking the $p$-th derivative of the expression above, we obtain
\begin{equation}
	y_p(t) = p! \sum_{n=0}^p   \binom{p}{p-n}\frac{a_n}{(N+2n)_{p-n}}e^{-\Lambda_n t}.
\end{equation}
The coefficients $a_n$ depend on the initial condition (see the previous section).
\end{proof}

From the explicit computations in the previous section, we can specify to the coefficients $a_n$ for the initial condition $\psi=(1,0\ldots,0)$. Indeed, the moments of $|U_t^{11}|^{2}$, are
\begin{equation}
	\E\left[|U_t^{11}|^{2p}\right] = p!\sum_{n=0}^p   \binom{p}{p-n}\frac{(N-1)_n }{(N+n-1)_n(N+2n)_{p-n}}e^{-\Lambda_n t}
\end{equation}
and the moments for $|U_t^{j1}|^{2}, j>1$, are given
\begin{equation}
	\E\left[|U_t^{j1}|^{2p}\right] = p!\sum_{n=0}^p   \binom{p}{p-n}\frac{(-1)^n }{(N+n-1)_n(N+2n)_{p-n}}e^{-\Lambda_n t}.
\end{equation}
For general initial conditions, one has to compute recursively all the moments $y_p$. In the next proposition we give the general form of $y_1$ and $y_2$ for general initial conditions.

\begin{prop} The moments $y_1$ and $y_2$ are given by
\begin{eqnarray}
y_1(t)&=&\left(y_1(0)-\frac 1N\right)e^{-t}+\frac 1N,\\
y_2(t)&=&\left[y_2(0)-\frac{1}{N+2}\left(y_1(0)-\frac 1N\right)-\frac {2}{N(N+1)}\right]e^{-(2+2/N)t}\nonumber\\&&+\frac{1}{N+2}\left[y_1(0)-\frac 1N\right]e^{-t}+\frac {2}{N(N+1)},
\end{eqnarray}
where $y_1(0)=|U_0^{j1}|^{2}=|\psi^{j}|^{2}$ and $y_2(0)=|U_0^{j1}|^{4}=|\psi^{j}|^{4}.$
\end{prop}

\bigskip

We now move on to study the covariance of elements $|U_t^{j1}|^{2}$. We aim to compute the quantities
 $$\E\left[|U_t^{j1}|^{2}|U_t^{k1}|^{2}\right],\,\, j \neq k=1,\ldots,N.$$
 Note that for $j=k$, one can use the more general moment formula in Proposition \ref{prop:moments}. Let us define define a family of covariance maps $f_{jk}:[0, \iy) \to [0, 1]$, with $f_{jk}(t) = \E\left[|U_t^{j1}|^{2}|U_t^{k1}|^{2}\right]$. Using the Laplace transform, we have
 \begin{equation}
f_{jk}(t)=\frac{\partial^2 \phi}{\partial\lambda_j\lambda_k}(0;t).
 \end{equation}
 We put also $f_j(t)=E \left[|U_t^{j1}|^{2}\right], \forall t\geq0$. From Proposition \ref{prop:moments}, we have that 
 $$f_j(t)=\partial_{\lambda_j}\phi(0;t)=\E \left[|U_t^{j1}|^{2}\right]=\left(f_j(0)-\frac 1N\right)e^{-t}+\frac 1N.$$  
 
 \begin{prop}\label{prop:covariances} 
 For all $j \neq k$, the covariance applications satisfy the following first order ODE
 \begin{equation}\label{eq-covar} 
 f'_{jk}=\frac{1}{N}\Big(f_j+f_k\Big)-\left(2+\frac{2}{N}\right)f_{jk},
 \end{equation}
which admits the solution
 \begin{eqnarray}
 f_{jk}(t)&=&\left[f_{jk}(0)-\frac{1}{N+2}\left(f_j(0)+f_k(0)-\frac{2}{N}\right)-\frac{1}{N(N+1)}\right]e^{-t(2+2/N)}\nonumber\\
 &&+\frac{1}{N+2}\left[f_j(0)+f_k(0)-\frac{2}{N}\right]e^{-t}+\frac{1}{N(N+1)},
 \end{eqnarray}
 with $f_j(0) = |\psi^j|^2$, $f_k(0) = |\psi^k|^2$ and $f_{jk}(0) = |\psi^j|^2|\psi^k|^2$.
 \end{prop}

\begin{proof}
To compute the two times derivative it is sufficient to consider the functions $$\phi(\lambda_j,\lambda_k;t)=\phi(0,\ldots,0,\lambda_j,0,\ldots,0,\lambda_k,0,\ldots,0;t),$$
with $\lambda_j$ (resp.~$\lambda_k$) in the $j$-th (resp.~$k$-th) position.
 Indeed, we have $f_{jk}(t)=\partial_{\lambda_j\lambda_k}\phi(0,0;t)$. Using the Laplace transform of Theorem \ref{thm:main}, we obtain
 \begin{eqnarray}
 \partial_t\phi&=&\frac{\lambda_j}{N}\phi+\frac{\lambda_k}{N}\phi+\left(\frac{\lambda_j^2}{N}-\lambda_j\right)\partial_{\lambda_j}\phi+\left(\frac{\lambda_k^2}{N}-\lambda_k\right)\partial_{\lambda_k}\phi\nonumber\\
 &&-\frac{1}{N}\Big(\lambda_j^2\partial_{\lambda_j\lambda_j}\phi+\lambda_k^2\partial_{\lambda_k\lambda_k}\phi+2\lambda_j\lambda_k\partial_{\lambda_j\lambda_k}\phi\Big).
 \end{eqnarray}
Applying $\partial_{\lambda_j\lambda_k}$ on both side and taking next $\lambda_j=\lambda_k=0$, we get the   expression
\begin{eqnarray}
f'_{jk}(t)=\frac{1}{N}\Big(\partial_{\lambda_k}\phi(0,0;t)+\partial_{\lambda_j}\phi(0,0;t)\Big)-\left(2+\frac{2}{N}\right)f_{jk}(t),
\end{eqnarray}
which corresponds exactly to the expression \eqref{eq-covar}.
\end{proof}

Naturally this expression depends on the initial condition. For example, in the case of $\psi=e_1$, we get
\begin{eqnarray}
f_{1j}(t)&=&\left(-\frac{1}{N+2}\left(1-\frac{2}{N}\right)-\frac{1}{N(N+1)}\right)e^{-t(2+2/N)}\nonumber\\
 &&+\frac{1}{N+2}\left(1-\frac{2}{N}\right)e^{-t}+\frac{1}{N(N+1)},\quad\mathrm{for}\,\,j\neq1\nonumber\\
f_{jk}(t)&=&\left(\frac{2}{(N+2)N}-\frac{1}{N(N+1)}\right)e^{-t(2+2/N)}\nonumber\\
 &&-\frac{2}{(N+2)N}e^{-t}+\frac{1}{N(N+1)}, \quad\mathrm{for}\,\,j\neq k\,\,\mathrm{and}\,\,j,k\neq1.
\end{eqnarray}

Finally, we estimate different statistics of interest in Quantum Information Theory. More precisely, we compute the average value of an observable and we give a bound on the average R\'enyi entropy of a pure with distribution $\pi^\psi_t$.

As stated in the introduction, the main motivation for this work was to define a new ensemble of random pure states. If a quantum system is in a state described by the vector $\psi_t$ having distribution  $\pi^\psi_t$ and an observable $A \in \M_N(\C)$ is measured, quantum theory predicts that the average value observed is $\langle \psi_t,A\psi_t\rangle$. 

We shall compute the average value $\mathbb{E}\langle \psi_t,A\psi_t\rangle = \mathbb{E}\langle \psi,U_t^*AU_t\psi\rangle = \mathbb{E}\langle \psi,U_tAU_t^*\psi\rangle$. 
\begin{lem}\label{lem:average-observable}
The average value of the measure of a fixed observable $A \in \M_N(\C)$ on a quantum system described by the ensemble $\pi^\psi_t$ is
\begin{equation}\label{eq:average-observable}
\mathbb{E}\langle \psi_t,A\psi_t\rangle=\left(\langle \psi,A\psi\rangle-\frac{\trace(A)}{N}\right)\,e^{-t}+\frac{\trace(A)}{N}
\end{equation}
\end{lem}
\begin{proof}
The result follows from the computation of $dU_tAU_t^*$ with the Ito calculus. Using the Ito formula \eqref{eq:ito}, we have
\begin{equation}
dU_tAU_t^*=U_t(d\,AU_t^*)+dU_t\,(AU_t^*)+\frac{\trace(A)}{N}I\d t.
\end{equation}
Computing $d\langle \psi,U_tAU_t^*\psi\rangle$ and taking the expectation, we get
\begin{equation}
d\mathbb{E} \langle \psi,U_tAU_t^*\psi\rangle =\left(-\mathbb{E} \langle \psi,U_tAU_t^*\psi\rangle +\frac{\trace(A)}{N}\right)\d t,
\end{equation}
which implies equation \eqref{eq:average-observable}.
\end{proof}

In particular, when $t$ goes to infinity, we recover the usual result for the Haar measure $\pi_\iy$. It is easy to compute average value of an observable $A \in \M_N(\C)$ for the Fubini-Study ensemble:
\begin{align}
	\langle A \rangle = \int \braobsket{\psi}{A}{\psi} \d\pi_\iy(\psi) &= \int \braobsket{e_1}{U^*AU}{e_1} \d{\Haar(U)}\\
	&= \braobsket{e_1}{\int U^*AU \d{\Haar(U)}}{e_1} = \\
	&= \braobsket{e_1}{\frac{\trace A}{N} I_N}{e_1} = \frac{\trace A}{N}.
\end{align}

We finish by deriving bounds for the R\'eyni entropies. Recall that these quantities are defined for a pure state $\psi$ by 
\begin{equation}
S_p(\psi)=\frac{1}{1-p}\ln\left[\sum_{j=1}^N|\psi^j|^{2p}\right],
\end{equation}
where $p\geq 2$ is an integer.
For the random pure state generated by a unitary Brownian motion, we have
\begin{equation}
S_p(\psi_t)=\frac{1}{1-p}\ln\left[\sum_{j=1}^N|U_t^{j1}|^{2p}\right]
\end{equation}
We shall estimate $\mathbb{E}[S_p(\psi_t)]$ with the help of the Jensen inequality. We have indeed
\begin{eqnarray}
\mathbb{E}[S_p(\psi_t)]&\geq&\frac{1}{1-p}\ln\left[\sum_{j=1}^N\mathbb{E}| U_t^{j1}|^{2p}\right] = \frac{1}{1-p}\ln Y_p,
\end{eqnarray}
where $Y_p(t)=\sum_{j=1}^N\mathbb{E}| U_t^{j1}|^{2p}$ are the sum-of-moments functions. Using the moment formulas obtained in the beginning of this section, one can compute in principle the function $Y_p(t)$. In particular, in the case when $\psi = e_1$, we have
\begin{equation}
	Y_p(t) = p!\sum_{n=0}^p   \binom{p}{p-n}\frac{(N-1)_n+(N-1)(-1)^n }{(N+n-1)_n(N+2n)_{p-n}}e^{-\Lambda_n t}.
\end{equation}

In the limit $t \to \infty$, only the term $n=0$ survives, and we obtain
$$\lim_{t \to \infty} Y_p(t) = (N+p){\binom{N+p}{p}}^{-1},$$
which is formula (7.60) in \cite{BZ06}. We use this expression to bound the entropy of the coordinate vector of an uniform point on the unit sphere of $\mathbb C^N$:
$$\mathbb E_{\pi_\infty}[S_p(\psi)] \geq \frac{\log \left[ \frac{1}{N+p}\binom{N+p}{p}\right]}{p-1}.$$

Taking the formal limit $p \to 1$ in the above expression, we obtain the known bound \cite[equation 7.70]{BZ06}
$$\mathbb E_{\pi_\infty}[S_p(\psi)] \geq \sum_{k=2}^N \frac{1}{k}.$$

\section*{Acknowledgments}

Both researchers were supported by a PEPS grant from the CNRS. I.~N.~ acknowledges financial support from the ANR project OSvsQPI 2011 BS01 008 01. C.~P.~ acknowledges financial support from the ANR project HAM-MARK, N${}^\circ$ ANR-09-BLAN-0098-01. We would also like to thank Myl\`ene Ma\"ida for interesting discussion around the unitary Brownian motion.

\end{document}